\documentclass{amsart}
\usepackage{amsmath,amssymb}
\newcommand{\Gal}{{\rm Gal}}
\newcommand{\gl}{{\rm GL}}
\newcommand{\pgl}{{\rm PGL}}
\renewcommand{\O}{\mathfrak O}
\newcommand{\p}{\mathfrak p}
\renewcommand{\P}{\mathfrak P}
\newcommand{\Q}{\mathbb Q}
\newcommand{\Z}{\mathbb Z}
\newcommand{\C}{\mathbb C}
\newcommand{\Vig}{Vign\'eras}
\theoremstyle{plain}
\newtheorem{thm}{Theorem}
\title{Even Icosahedral Galois Representations of Prime Conductor}
\author{Darrin Doud}
\address{Department of Mathematics, Brigham Young University, Provo, UT  84601}
\email{doud@math.byu.edu}
\author{Michael W. Moore}
\address{Department of Mathematics, Brigham Young University, Provo, UT  84601}
\email{mike@math.byu.edu}
\date{\today}
\begin{document}
\begin{abstract}
In this paper, we use a series of targeted Hunter searches to prove that the minimal prime conductor of an even icosahedral Galois representation is 1951.  In addition, we give a complete list of all even icosahedral Galois representations of prime conductor less than 10,000.
\end{abstract}

\maketitle

\section{Introduction}
Galois representations and their connections to modular forms have played a key role in much of modern number theory.  Odd two-dimensional Galois representations and their connection to holomorphic modular forms were important in Wiles' proof of Fermat's Last Theorem. Even two-dimensional Galois representations are conjectured to correspond to Maass wave forms.  Among the simplest examples of Galois representations are those which are ramified at only one prime.  In the two-dimensional case, they have been studied by Serre~\cite{Serre} and \Vig~\cite{Vigneras}.  Serre studied odd irreducible Galois representations with prime conductor, and determined that, besides the dihedral ones, there were three major classes of such representations.  \Vig\  extended Serre's study to even Galois representations and classified them similarly.  She gave several examples of dihedral and octahedral representations but was unable to find any even icosahedral Galois representation of prime conductor.  This is no surprise, as the discriminants of the number fields describing such  representations are far larger than $2\times 10^7$, and hence cannot be found even in modern tables of number fields~\cite{tables}.  
 In this paper, we describe the techniques that we used to find the even icosahedral representation of minimal prime conductor.  An analogous result for odd icosahedral Galois representations (namely that the minimal prime conductor is 2083) can be found in~\cite{Basmaji}. We also give a table of all even icosahedral representations of prime conductor having conductor less than 10,000.

\section{Even Galois representations}
A Galois representation is a continuous homomorphism $\rho:G_\Q\to \gl_2(\C)$. Continuity in this case simply means that the image of $\rho$ is finite.  Given a Galois representation $\rho$, we may construct a projective representation $\tilde\rho=\pi\circ\rho$, where $\pi:\gl_2(\C)\to\pgl_2(\C)$ is the natural projection map.  Such a projective representation has as its image a finite subgroup of $\pgl_2(\C)$.  Any such subgroup is well known to be cyclic, dihedral, $A_4$ (tetrahedral), $S_4$ (octahedral), or $A_5$ (icosahedral).  

A theorem of Tate~\cite{Serre} shows that given any continuous homomorphism $\tilde\rho:G_\Q\to\pgl_2(\C)$, there is a Galois representation $\rho:G_\Q\to\gl_2(\C)$ such that $\tilde\rho=\pi\circ\rho$. In fact, this lift is not uniquely defined but can be chosen to have almost any desired ramification properties.  In particular, $\rho$ can be chosen to be ramified only at those primes at which $\tilde\rho$ is ramified.  

Complex conjugation is a well defined element of $\Gal(\overline{\Q}/\Q)$, which we will denote by $\tau$.  A Galois representation is called odd if $\det(\rho(\tau))=-1$; it is called even if $\det(\rho(\tau))=1$.

\Vig~\cite{Vigneras} classified even Galois representations of prime conductor.  Her description involved the fixed field $E$ of the kernel of $\tilde\rho$ and the ramification index $e_p$ of $p$ in $E/\Q$.  She showed that the irreducible even Galois representations $\rho$ with prime conductor $p$ are either dihedral or correspond to field extensions $E/\Q$ satisfying

\begin{list}{}{}
\item[1.] $\Gal(E/\Q)\cong A_4$, $p\equiv 1\pmod 3$, $E/\Q$ is the normal closure of a quartic field of discriminant $p^2$, and $e_p=3$,
\item[2a.] $\Gal(E/\Q)\cong S_4$, $p\equiv 5\pmod 8$, $E/\Q$ is the normal closure of a quartic field of discriminant $p$, and $e_p=2$,
\item[2b.] $\Gal(E/\Q)\cong S_4$, $p\equiv 1\pmod 8$, $E/\Q$ is the normal closure of a quartic field of discriminant $p$ or $p^3$, and $e_p=2$ or $e_p=4$,
\item[3a.] $\Gal(E/\Q)\cong A_5$, $p\equiv 1\pmod 5$,  $E/\Q$ is the normal closure of a quintic field of discriminant $p^4$, and $e_p=5$,
\item[3b.] $\Gal(E/\Q)\cong A_5$, $p\equiv 1\pmod 3$, $E/\Q$ is the normal closure of a quintic field of discriminant $p^2$, and $e_p=3$,
\item[3c.] $\Gal(E/\Q)\cong A_5$, $p\equiv 1\pmod 4$, $E/\Q$ is the normal closure of a quintic field of discriminant $p^2$, and $e_p=2$.
\end{list}

Note that determining a Galois extension $E$ with the properties described here is equivalent to defining a projective representation $\tilde\rho$, so that finding even Galois representations of prime conductor is reduced to finding number fields with certain properties.

\Vig\  indicated that there were no even representations of prime conductor of type 1 with $p<106$; in fact, checking  tables of quartic number fields~\cite{tables} (which were not available to \Vig) we find that the smallest $p$ for which such a representation exists is $p=163$.

\Vig\  gave several examples of representations of type $2$.  Using a combination of class field theory and number field tables \cite{tables}, we may determine that the minimal prime conductor of a representation of type 2a is 8069, and of a representation of type 2b is 2713.

\Vig\  gave no  examples of representations of type 3.  She indicates that the tables of Buhler~\cite{Buhler} contain no such representations.  These tables, which contain representations of conductor up to 10,000, make no claim of completeness, so that the possibility remains that such representations exist.
Modern tables of number fields \cite{tables}, which give all totally real quintic number fields of discriminant at most $2\times 10^7$ yield no examples of such a representation.  This shows that any representation of type 3b or 3c must have $p>4472$ (since $p^2>2\times 10^7$), and that any representation of type 3a must have $p>66$ (since $p^4>2\times 10^7$).  

Representations of type 3a with conductor much smaller than 10,000, but corresponding to number fields of discriminant much greater than $2\times10^7$ do exist.  Because the discriminants of the fields are so large, the fields cannot be found in any standard tables of number fields.  We will use a series of targeted Hunter searches to find all representations of type 3a with prime conductor less than 10,000. In addition, we will show that there are no representations of types 3b and 3c with conductor less than 10,000 and determine the minimal conductor of each of the three types of representation.

\section{Targeted Hunter Searches}

Given a discriminant $D$ and a degree $n$, the goal of a Hunter search is to find a polynomial defining a number field $K$ of degree $n$ and discriminant $D$ or to determine that no such number field exists.  Such a search is conducted by using Hunter's theorem and various relations between coefficients of polynomials to bound the coefficients of a polynomial defining a number field with the given degree and discriminant.  Such a search generally yields a set of possible polynomials which is far too large to examine in a reasonable time (especially for large $n$ or $D$).  A refinement of the Hunter search can be used if we are only interested in searching for $K$ with a desired ramification type at a certain prime $p$.  We use our knowledge of the ramification of the desired field $K$ to determine congruence conditions on the defining polynomial of $K$.  This allows us to drastically decrease the number of polynomials which need to be investigated.  This ``targeting'' of a specific ramification structure is a powerful tool and can make otherwise impractical searches quite feasible. We describe both the traditional Hunter search and also the details behind a targeted search in this paper.  These techniques are not original; they were described in~\cite{Cohen2} and~\cite{JR} in some detail, and our purpose in describing them here is to use them to find even icosahedral representations of prime conductor.

\subsection{Hunter's Theorem}
Hunter's theorem~\cite[Theorem 6.4.2]{Cohen1} states:
\begin{thm}{}
Let $K$ be a number field of degree $n$ over $\Q$.  There exists $\theta\in\O_K\setminus\Z$ such that
$$\sum_{i=1}^n|\theta_i|^2\leq \frac1n{\rm{Tr}}(\theta)^2+\gamma_{n-1}\left(\frac{|d(K)|}n\right)^{1/(n-1)}$$
where $d(K)$ is the discriminant of $K$, ${\rm Tr}(\theta)$ is the trace of $\theta$, $\gamma_{n-1}$ is Hermite's constant~\cite[p. 20]{Conway} in dimension $n-1$, and $\theta_1,\ldots,\theta_n$ are the conjugates of $\theta$.  In addition, we may assume that $0\leq{\rm{Tr}}(\theta)\leq n/2$.
\end{thm}

If $K$ is a primitive field, then the $\theta\in K$ which is guaranteed, by Hunter's theorem, to exist is in fact a generator of $K$ over $\Q$; i.e. $K=\Q(\theta)$.  Then the minimal polynomial of $\theta$ is a defining polynomial for $K$.  

In order to conduct a Hunter search for a primitive number field with a given discriminant $D$ and degree $n=5$, we first assume that such a number field exists.  By Hunter's theorem there will then be a $\theta\in K$ with conjugates $\theta=\theta_1,\ldots,\theta_5$, such that
$$\sum_{i=1}^5|\theta_i|^2\leq \frac15{\rm{Tr}}(\theta)^2+\gamma_{4}\left(\frac{|D|}5\right)^{1/4}=t_2$$
and $0\leq{\rm Tr}(\theta)\leq 5/2$. We now write the minimal polynomial $f$ of $\theta$ as
$$f(x)=x^5-a_1x^{4}+a_2x^{3}-\cdots-a_5,$$
with $a_i\in\Z$ for $1\leq i\leq 5$, and we see that $a_1= {\rm Tr}(\theta)$, so $0\leq a_1\leq 5/2$. For future reference, we will define $a_0=1$.

For $k>0$, we now define 
$$S_k=\sum_{i=1}^5\theta_i^k\quad{\text{and}}\quad T_k=\sum_{i=1}^5|\theta_i|^k.$$  We note that the $S_k$ are symmetric polynomials in the roots of $f$, and are thus integers.  Now for $k>1$, we see~\cite{Cohen2} that $|S_k|\leq T_k<T_2^{k/2}<t_2^{k/2}$.  We will also use Newton's relations, which state that
$$ka_k=\sum_{j=1}^k(-1)^{j-1}a_{k-j}S_j.$$
In the case $k=2$, solving Newton's relations for $a_2$  shows that
$$a_2=\frac{\left(a_1S_1-a_0S_2\right)}2,$$
and since $-t_2\leq S_2\leq t_2$ and $S_1=a_1$, we see that
$$\frac{a_1^2-t_2}2\leq a_2\leq\frac{a_1^2+t_2}2.$$
We may improve this bound by using the fact that the field $K$ for which we are searching must be totally real.  In this case,~\cite[Prop. 9.3.7(2)]{Cohen2} shows that $a_2\leq(a_1^2-6)/2$. Hence, each possible value for $a_1$ gives a finite range of possibilities for $a_2$.

Similarly, solving Newton's relation with $k=3$ for $a_3$, we obtain
$$a_3=\frac{a_2a_1-a_1S_2+a_0S_3}{3}.$$
Now, note that $|S_3|\leq T_3\leq S_2^{3/2}$.  In addition, we know that $S_2=a_1^2-2a_2$.  Putting these into the expression for $a_3$, we see that
$$\frac{-a_1^3+3a_1a_2-S_2^{3/2}}{3}\leq a_3\leq\frac{-a_1^3+3a_1a_2+S_2^{3/2}}{3},$$
so that each possible pair of values for $a_1$ and $a_2$ gives a finite range of values for $a_3$.

At this point, it is convenient to bound $a_5$ before bounding $a_4$.  By~\cite[Prop. 9.3.7(1)]{Cohen2} we see that $|a_5|< S_2^{5/2}/5^{5/2}$, and note that this bound depends only on the values of $a_1$ and $a_2$.  Finally, solving Newton's relation with $k=4$ for $a_4$, we obtain
$$a_4=\frac{a_1(a_3+S_3)-a_2S_2-S_4}{4}.$$
Using the fact that $5|a_5|^{4/5}\leq S_4\leq S_2^2$~\cite[Prop. 9.3.7(3)]{Cohen2}, we obtain the inequalities
$$\displaystyle{\frac{(a_3+S_3)a_1-a_2S_2-S_2^2}4}\leq a_4\leq\displaystyle{\frac{(a_3+S_3)a_1-a_2S_2-5|a_5|^{4/5}}4}.$$

At this point, since each $a_k$ is constrained to lie in a finite interval, we have constructed a finite set of polynomials, which contains a defining polynomial for $K$ (under the assumption that $K$ exists).  We now need only check each polynomial, and find the ones defining a number field of the correct discriminant.  If no such polynomial exists in the search space, then we are guaranteed that no totally real field with the desired degree and discriminant can exist.

Although our search space is finite, for even moderately large $D$ these bounds give a search space which is far too large to search exhaustively. It is possible to get better bounds on the coefficients by more sophisticated methods~\cite{Pohst}, but we chose instead to reduce our search space by using ramification properties of the desired field $K$.

\subsection{Targeted searches}
Hunter's theorem gives us a finite (although very large) search space which is guaranteed to contain a polynomial defining the field for which we are searching (if such a field exists).  Targeting is a technique to take advantage of the fact that we may know something about the ramification of the fields for which we are searching.

\begin{thm} 
Let $L/K$ be a finite extension of number fields such that $L=K(\alpha)$, with $\alpha$ a root of a monic irreducible polynomial $f(x)\in \O_K[x]$.  Let $\p$ be a prime ideal of $\O_K$, and suppose that $\p\O_L=\prod\P_i^{e_i}$ is a product of powers of distinct prime ideals in $L$, such that the inertial degree of $\P_i$ is $f_i$.  Then $f\equiv\prod g_i^{e_i}\pmod\p$ with $g_i\in \O_K[x]$ and $\deg(g_i)=f_i$.
\end{thm}

{\em Remark:} Note that there is no guarantee that the $g_i$ are irreducible. In order to guarantee their irreducibility, we would need to know something about the index of $\O_K[\alpha]$ in $\O_L$.

\begin{proof} Denote by $K_\p$ and $L_{\P_i}$ the completions of $K$ and $L$ with respect to the indicated primes.  We know~\cite{CF} that 
$$K_\p\otimes_K L\cong\bigoplus_{\P_i}L_{\P_i},$$
which immediately implies that in $K_\p$, we have $f=\prod h_i$, with each $h_i$ having degree $e_if_i$, and $L_{\P_i}\cong K_\p[x]/(h_i(x))$~\cite[page 58]{CF}.  We thus only need to show that each $h_i$ is congruent modulo $\p$ to some $g_i^{e_i}$.

Let $M$ be a splitting field of $h_i$ over $K_\p$, let $\P$ be the unique prime ideal of $M$, let $I$ be the inertia group of $M/K_\p$, and let $\beta$ be any root of $h_i$ in $M$.  Let $e$ and $f$ be the ramification index and inertial degree of $M/K_\p$, and let $e_i$ and $f_i$ be the ramification index and inertial degree of $K_\p(\beta)/K_\p$.  Note that the values of $e_i$ and $f_i$ are the same as the $e_i$ and $f_i$ defined above, and that they are independent of the choice of $\beta$.   Let $I_\beta$ be the inertia group of $M/K_\p(\beta)$.  Then $I_\beta$ has order $e'=e/e_i$, and is equal to $I\cap\Gal(M/K_\p(\beta))$. If we consider the action of $I$ on the roots of $h_i$, we see that the stabilizer of $\beta$ is just $I_\beta$.  Then the orbit-stabilizer relationship shows that the orbit of $\beta$ under the action of $I$ has order $e_i$.  However, the image of $\beta$ under any element of $I$ is congruent to $\beta$ modulo $\P$. Hence each root $\beta$ of $h_i$ is congruent modulo $\P$ to the $e_i$ roots in its orbit.  Let $\beta_1,\ldots,\beta_{f_i}$ be a set of representatives of the orbits.  Define
$$p(x)=\prod_{j=1}^{f_i}(x-\beta_j).$$
Clearly the coefficients of $p(x)$ are fixed modulo $\P$ by $\Gal(M/K_\p)$, so that the reduction $\bar p(x)$ modulo $\P$ actually has coefficients in $\O_{K_\p}/\p$.  Hence, $\bar p(x)$ is actually the reduction modulo $\p$ of a polynomial $g_i(x)\in\O_K[x]$, and we see that $h_i(x)\equiv g_i(x)^{e_i}\pmod p$, as desired.\end{proof}

For each of the three types of even icosahedral Galois representation described by \Vig, we need to search for a totally real quintic field $K$ ramified only at one prime $p$, with a specific type of ramification at $p$.  We do this by conducting a Hunter search for a polynomial $f$ which defines $K$.  Since we know the desired discriminant and signature of $K$, we know that a defining polynomial exists in a finite (but large) set of polynomials.  In addition, because we know the desired ramification of $p$ in $K$, we know how this defining polynomial must factor modulo $p$.  We give details for each type of representation.

{\bf Type 3a}: Here $p\equiv 1\pmod 5$ is ramified with $e=5$.  Hence, the desired polynomial is of the form $(x-a)^5\pmod p$.  In our program, as soon as we determine $a_1$, we know exactly what $a$ is (since $a_1\equiv 5a\pmod p$), and this determines the values of $a_2,\ldots,a_5$ modulo $p$.  Rather than examining all values of $a_2,\ldots,a_5$ within our bounds, we only examine those satisfying the given congruence conditions.  This results in examining approximately 1 of every $p^4$ polynomials in the search space, dramatically reducing the time required for the search.

{\bf Type 3b}: Here $p\equiv 1\pmod 3$ and $p$ is ramified with $e=3$.  There may be either one or two other primes above $p$ (besides the ramified one), with inertial degrees either $2$ or $1$.  In any case, we see easily that $f\equiv (x-a)^3(x^2+bx+c)\pmod p$, where the quadratic factor may be either reducible or irreducible.  A we proceed through our Hunter search then, we determine $a_1,a_2,$ and $a_3$.  This allows us to solve for $a$ (as a root of a cubic polynomial modulo $p$) and afterward determine $b$ and $c$.  Once $a$, $b$, and $c$ are determined, the values of $a_4$ and $a_5$ are determined modulo $p$.  Hence, we only need to examine about 1 of every $p^2$ polynomials in the search space.

{\bf Type 3c}: Here $p\equiv 1\pmod 4$ and $p$ splits in $K$ as either $P^2Q$ with $P$ having inertial degree $2$ or as $P^2Q^2R$ with all primes having inertial degree 1.  In the first case, $f$ factors modulo $p$ as $(x^2+ax+b)^2(x+c)$, and in the second case, $f$ factors modulo $p$ as $(x+d)^2(x+e)^2(x+c)$.  In either case we may write $f\equiv(x^2+ax+b)^2(x+c)\pmod p$, where the quadratic factor may be reducible.  As above, determining the first three coefficients of $f$ determines $a,b,$ and $c$ modulo $p$, which determine $a_4$ and $a_5$ modulo $p$.  Hence, we need only examine $1$ of every $p^2$ polynomials in our search space.

\section{Computational Results}

We programmed the targeted hunter searches described above in the GP/Pari \cite{Pari} programming language.  For each prime congruent to 1 modulo 5 and below 10,000 we searched for an $A_5$-extension corresponding to a representation of type 3a.  The results of our search are summarized in table \ref{table}.  
\begin{table}
\begin{tabular}{|c|c|}
\hline
Conductor&Defining Polynomial for $E/\Q$\cr
\hline
1951&$x^5 - x^4 - 780x^3 - 1795x^2 + 3106x + 344$\cr
\hline
2141&$x^5 - x^4 - 856x^3 + 4025x^2 + 28501x - 40877$\cr
\hline
3701&$x^5 - x^4 - 1480x^3 - 18209x^2 + 2191x + 9683$\cr
\hline
3821&$x^5 - x^4 - 1528x^3 - 1987x^2 + 16629x - 12281$\cr
\hline
8501&$x^5 - x^4 - 3400x^3 - 41825x^2 + 671511x - 966731$\cr
\hline
9461&$x^5 - x^4 - 3784x^3 + 2649x^2 + 2960082x - 2781864$\cr
\hline
\end{tabular}
\caption{Defining polynomials for all even icosahedral Galois representations of prime conductor less than $10,000$.}\label{table}
\end{table}
We stress that for primes $p$ less than 10,000 and not included in this table, no even icosahedral Galois representation of conductor $p$ can exist.  We found a total of six primes $p$ for which there exists an even icosahedral Galois representation of conductor $p$.  For each of these primes we list a defining polynomial for the quintic field $K$.  We note that of these extensions, at least two were previously known~\cite[Section 5.1.3]{ADP}; however, the extension with minimal conductor seems not to have been known.   From the table, we see that the minimal prime conductor of an even icosahedral Galois representation is 1951.  

This sequence of searches took approximately 5
processor-weeks to complete on a dual processor Xeon 1.7 GHz
machine. Without the improvement of
targeting, the computational time required would have made such a
calculation unrealistic even if we had used a much faster machine.

We also conducted a sequence of searches for representations of type 3b and 3c, using the targeting techniques described above.  We found that the smallest prime $p$ for which there is a representation of type 3b is $p=\text{10,267}$, and the smallest prime $p$ for which there is a representation of type 3c is $p=\text{13,613}$. Each of these sequences of searches took about two to three hours.  This time is dramatically shorter than the searches for type 3a representations, primarily because of the smaller size of the discriminants of the fields involved.

As a check that the program indeed searches all polynomials in the search space, we note that for every prime $p$ congruent to 1 modulo 5 there is a $\Z/5\Z$-extension of $\Q$ ramified only at $p$ and having discriminant $p^4$, namely the degree 5 subfield of the $p$th cyclotomic field.  By Hunter's theorem, our search space should contain a polynomial defining this field.  In fact for every such prime $p$ within our search range, we did find a defining polynomial for this field. 

A second check is that for every prime congruent to 1 modulo 4 such that the class number of $\Q(\sqrt p)$ is divisible by 5 there is a dihedral extension of degree 10 ramified only at $p$.  Since this dihedral extension will be defined by a quintic polynomial satisfying the same ramification and discriminant conditions as  the quintic polynomial for a type 3c representation, it should also be found by our program.  In fact, for every such prime within our search range we did find a defining polynomial for this dihedral extension.

\bibliographystyle{amsplain}

\providecommand{\bysame}{\leavevmode\hbox to3em{\hrulefill}\thinspace}

\end{document}